\documentclass[a4paper]{amsart}

\usepackage{amsmath,amsthm,amsfonts,amssymb,amscd}
\input{xy}
\xyoption{all}

\newtheorem{theo}{Theorem}[section]
\newtheorem{pro}[theo]{Proposition}

\newtheorem{defi}[theo]{Definition}

\newtheorem*{exa}{Example}

\newtheorem{prob}{Problem}

\newtheorem*{ackno}{Acknowledgements}

\newenvironment{sis}{\left\{\begin{aligned}}{\end{aligned}\right.}

\numberwithin{equation}{section}

\newcommand{\N}{\mathbb{N}}

\renewcommand{\d}{{\rm d}}

\newcommand{\ti}{\widetilde}

\newcommand{\ep}{\epsilon}

\newcommand{\Sq}{{\rm Sq}}

\begin{document}

\title{Restricted simple Lie algebras and their infinitesimal deformations}
\author{Filippo Viviani}
\address{Universit\'a degli studi di Roma Tor Vergata, Dipartimento di
Matematica, via della Ricerca Scientifica 1, 00133 Rome}
\email{viviani@mat.uniroma2.it}

\thanks{The author was supported by a grant from the Mittag-Leffler Institute of Stockholm.}


\begin{abstract}
In the first two sections, we  
review the Block-Wilson-Premet-Strade classification of restricted simple Lie algebras.
In the third section, we compute their infinitesimal deformations.
In the last section, we indicate some possible generalizations by formulating some open 
problems.
\end{abstract}

\maketitle

\section{Restricted Lie algebras}

We fix a field $F$ of characteristic $p>0$ and we denote with $\mathbb{F}_p$ the prime field
with $p$ elements.  
All the Lie algebras that we will consider are of finite dimension over $F$.
We are interested in particular class of Lie algebras, called restricted 
(or $p$-Lie algebras). 

\begin{defi}[Jacobson \cite{JAC2}]
A Lie algebra $L$ over $F$ is said to be restricted (or a $p$-Lie algebra) if there exits
a map (called $p$-map), $[p]:L\to L$, $x\mapsto x^{[p]}$, 
which verifies the following conditions:
\begin{enumerate}
\item ${\rm ad}(x^{[p]})={\rm ad}(x)^{[p]}$ for every $x\in L$.
\item $(\alpha x)^{[p]}=\alpha^p x^{[p]}$ for every $x\in L$ and every $\alpha\in F$.
\item $(x_0+x_1)^{[p]}=x_0^{[p]}+x_1^{[p]}+\sum_{i=1}^{p-1} 
s_i(x_0,x_1)$ for every $x,y\in L$,
where the element $s_i(x_0,x_1)\in L$ is defined by
$$s_i(x_0,x_1)=-\frac{1}{r}\sum_u {\rm ad}x_{u(1)}\circ {\rm ad}x_{u(2)}\circ
\cdots \circ {\rm ad}x_{u(p-1)}(x_1),$$
the summation being over all the maps $u:[1,\cdots,p-1]\to \{0,1\}$ taking $r$-times
the value $0$.
\end{enumerate}
\end{defi}

\begin{exa}
\begin{enumerate}

\item Let $A$ an associative $F$-algebra. Then the Lie algebra ${\rm Der}_F A$ of 
$F$-derivations of $A$ is a restricted Lie algebra with respect to the $p$-map
$D\mapsto D^p:=D\circ \cdots\circ D$.

\item Let $G$ a group scheme over $F$. Then the Lie algebra ${\rm Lie}(G)$ associated 
to $G$ is a restricted Lie algebra with respect to the $p$-map given by the differential 
of the homomorphism $G\to G$, $x\mapsto x^p:=x\circ \cdots\circ x$.

\end{enumerate}
\end{exa}

One can naturally ask when a $F$-Lie algebra can acquire the structure of a restricted
Lie algebra and how many such structures there can be. The following criterion 
of Jacobson answers to that question.

\begin{pro}[Jacobson]
Let $L$ be a Lie algebra over $F$. Then
\begin{enumerate}
\item It is possible to define a $p$-map on $L$ 
if and only if, for every element $x\in L$, the $p$-th iterate of ${\rm ad}(x)$ is
still an inner derivation. 
\item Two such $p$-maps differ by a semilinear map from $L$ to the center $Z(L)$ of $L$, 
that is a map  $f:L\to Z(L)$ such that
$f(\alpha x)=\alpha^p f(x)$ for every $x\in L$ and $\alpha\in F$.
\end{enumerate}
\end{pro}
\begin{proof}
See \cite[Chapter V.7]{JACLIE}.
\end{proof}

Many of the modular Lie algebras that arise ``in nature'' are restricted.
As an example of this principle, we would like to recall the following two results 
from the theory of finite group schemes and the theory of inseparable field extensions.

\begin{theo}\label{height-1}
There is a bijective correspondence 
$$\{\text{Restricted Lie algebras}/ F\} \longleftrightarrow 
\{\text{Finite group schemes}/F \text{ of height }1 \},$$
where a finite group scheme $G$ has height $1$ if the Frobenius $F:G\to G^{(p)}$ 
is zero. Explicitly to a finite group scheme $G$ of height $1$, one associates 
the restricted Lie algebra ${\rm Lie}(G):=T_0 G$. Conversely, to a restricted Lie
algebra $(L,[p])$, one associates the finite group scheme corresponding to the dual
of the restricted enveloping Hopf algebra $U^{[p]}(L):=U(L)/(x^p-x^{[p]})$. 
\end{theo}
\begin{proof}
See \cite[Chapter 2.7]{DG}.
\end{proof}

\begin{theo}
Suppose that $[F:F^p]<\infty$.
There is a bijective correspondence 
$$\{\text{Inseparable subextensions of exponent } 1\}  
\longleftrightarrow \{\text{Restricted subalgebras of }Der(F)\}$$ 
where the inseparable subextensions of exponent $1$ are the subfields
$E\subset F$ such that $F^p\subset E\subset F$ and ${\rm Der}(F):={\rm Der}_
{\mathbb{F}_p}(F)={\rm Der}_{F^p}(F)$. Explicitly to any field $F^p\subset E\subset F$
one associates the restricted subalgebra ${\rm Der}_E(F)$. Conversely, to any 
restricted subalgebra $L\subset {\rm Der}(F)$, one associates the subfield 
$E_L:=\{x\in F\: |\: D(x)=0 \text{ for all }D\in L\}$.
\end{theo}
\begin{proof}
See \cite[Chapter 8.16]{JACALG}. 
\end{proof}

\section{Classification of restricted simple Lie algebras}

Simple Lie algebras over an algebraically closed field of {\bf characteristic zero} were
classified at the beginning of the XIX century by Killing and Cartan.
The classification proceeds as follows:
first the non-degeneracy of the Killing form is used to establish a correspondence between 
simple Lie algebras and irreducible root systems and  then the irreducible root systems  are 
classified by mean of their associated Dynkin diagrams.  It turns out that 
there are four infinite families of Dynkin diagrams, called $A_n$, $B_n$, 
$C_n$, $D_n$, and five exceptional Dynkin diagram, called $E_6$, $E_7$, $E_8$, $F_4$ and 
$G_2$.
The four infinite families correspond, respectively, to the 
the special linear algebra $\mathfrak{sl}(n+1)$, the special orthogonal algebra of odd rank 
$\mathfrak{so}(2n+1)$, the symplectic algebra $\mathfrak{sp}(2n)$ and the special orthogonal
algebra of even rank $\mathfrak{so}(2n)$. For the simple Lie algebras corresponding to the 
exceptional Dynkin diagrams, see the book \cite{JACEXE} or the nice account in \cite{BAE}.

These simple Lie algebras admits a model over the integers via the (so-called) Chevalley bases.
Therefore, via reduction modulo a prime $p$, one obtains a restricted 
Lie algebra over $\mathbb{F}_p$, which is simple up to a quotient by a small ideal. 
For example $\mathfrak{sl}(n)$ is not simple 
if $p$ divide $n$, but its quotient $\mathfrak{psl}(n)=\mathfrak{sl}(n)/(I_n)$ by the unit 
matrix $I_n$ becomes simple.
There are similar phenomena occuring only for $p=2, 3$ for the other Lie 
algebras (see \cite[Page 209]{STR} or \cite{SEL}). 
The restricted simple algebras obtained in this way are 
called algebras of {\bf classical type}.  
Their Killing form is non-degenerate except at a finite number of primes.
Moreover, they can be characterized 
as those restricted simple Lie algebras admitting a projective representation 
with nondegenerate trace form (see \cite{BLO}, \cite{KAP}).

However, there are restricted simple Lie algebras which have
no analogous in characteristic zero and therefore are called nonclassical.
The first example of a nonclassical restricted simple Lie algebra is due to E. Witt, 
who in 1937 
realized that the derivation algebra $W(1):= {\rm Der}_F(F[X]/(X^p))$ over a field $F$ of 
characteristic $p>3$ is simple with a degenerate Killing form.
In the succeeding three decades, many more nonclassical restricted simple Lie algebras 
have been found (see \cite{JAC1}, \cite{FRA1}, \cite{AF}, \cite{FRA2}). 
The first comprehensive conceptual approach to constructing these 
nonclassical restricted simple Lie algebras was proposed by Kostrikin and Shafarevich in 1966
(see \cite{KS}). 
They showed that all the known examples  can be 
constructed as finite-dimensional analogues of the four 
classes of infinite-dimensional complex simple Lie algebras, which occurred in Cartan's 
classification of Lie pseudogroups (see \cite{CAR}).
These restricted simple Lie algebras, 
called of {\bf Cartan-type}, are divided into four families, called Witt-Jacobson, Special,
Hamiltonian and Contact algebras. 

\begin{defi}
Let $A(n):=F[x_1,\cdots,x_n]/(x_1^p, \cdots, x_n^p)$ the algebra of 
$p$-truncated polynomials in $n$ variables. Then the Witt-Jacobson Lie algebra $W(n)$
is the derivation algebra of $A(n)$:
$$W(n)={\rm Der}_F A(n).$$
\end{defi}

For every $j\in \{1, \dots , n\}$, we put $D_j:=\frac{\partial}{\partial x_j}$.
The Witt-Jacobson algebra $W(n)$ is
a free $A(n)$-module with basis $\{D_1, \dots, D_n\}$. Hence ${\rm dim}_F W(n) =np^n$ with a
basis over $F$ given by $\{x^aD_j\: | \: 1\leq j\leq n, \: x^a\in A(n)\}$.

The other three families are defined as $m$-th derived algebras of the subalgebras 
of derivations fixing a volume form, a Hamiltonian form and a contact form, respectively. 
More precisely, consider the natural action of $W(n)$ on the exterior algebra of differential
forms in $\d x_1,\cdots, \d x_n$ over $A(n)$. 
Define the following three forms, called volume form, Hamiltonian form and contact form:
$$\begin{sis}
& \omega_S=\d x_1\wedge \cdots \d x_n,\\
& \omega_H=\sum_{i=1}^m \d x_i \wedge \d x_{i+m} \: \text{ if } n=2m,\\
& \omega_K=\d x_{2m+1}+\sum_{i=1}^m (x_{i+m} \d x_i -x_i \d x_{i+m}) \: \text{ if }n=2m+1.
\end{sis}$$

\begin{defi}
Consider the following three subalgebras of $W(n)$:
$$\begin{sis}
& \widetilde{S}(n)=\{D\in W(n) \: |\: D\omega_S=0\},\\
& \widetilde{H}(n)=\{D\in W(n)\: | \: D \omega_H=0\},\\
& \widetilde{K}(n)=\{D\in W(n)\: | \: D \omega_K \in A(n)\omega_K\}.
\end{sis}$$
Then the Special algebra $S(n)$ ($n\geq 3$) is the derived algebra of 
$\widetilde{S}(n)$, while the Hamiltonian algebra $H(n)$
($n=2m\geq 2$) and the Contact algebra  $K(n)$ ($n=2m+1\geq 3$)
are the second derived algebras of $\widetilde{H}(n)$ and $\widetilde{K}(n)$, 
respectively.
\end{defi}

We want to describe more explicitly the above algebras, starting from 
the Special algebra $S(n)$. For every $1\leq i, j\leq n$ consider the following maps
$$D_{ij}=-D_{ji}:\begin{sis}
A(n)& \longrightarrow W(n)\\
f & \mapsto D_j(f)D_i-D_i(f)D_j.\\
\end{sis}$$

\begin{pro}
The algebra $S(n)$ has $F$-dimension  equal to $(n-1)(p^n-1)$ and 
is generated by the elements $D_{ij}(x^a)$ for $x^a\in A(n)$
and $1\leq i<j\leq n$.
\end{pro}
\begin{proof}
See \cite[Chapter 4.3]{FS}.
\end{proof}

Suppose now that $n=2m\geq 2$ and consider the map $D_H:A(n)\to W(n)$ 
defined by
$$D_H(f)=\sum_{i=1}^m\left[D_i(f)D_{i+m}-D_{i+m}(f)
D_i\right],$$
where, as before, $D_i:=\frac{\partial}{\partial x_i}\in W(n)$. Then the Hamiltonian
algebra can be described as follows:
\begin{pro}
The above map $D_H$ induces an isomorphism 
$$D_H: A(n)_{\neq 1, x^{\sigma}} \stackrel{\cong}{\longrightarrow} H(n),$$
where $A(n)_{\neq 1, x^{\sigma}}=\{x^a\in A(n)\: | \: x^a\neq 1, x^a\neq x^{\sigma}:=
x_1^{p-1}\cdots x_n^{p-1}\}.$ Therefore $H(n)$ has dimension $p^n-2$.
\end{pro}
\begin{proof}
See \cite[Chapter 4.4]{FS}.
\end{proof}

Suppose finally that $n=2m+1\geq 3$. Consider the map $D_K:A(n)\to K(n)$ defined by
$$D_K(f)=\sum_{i=1}^m\left[D_i(f)D_{i+m}-D_{i+m}(f)D_i\right]
+\sum_{j=1}^{2m} x_j\left[D_n(f)D_j-D_j(f)D_n\right] +2f D_n.
$$
Then the Contact algebra can be described as follows:

\begin{pro}
The above map $D_K$ induces an isomorphism
 $$K(n) \cong
\begin{sis}
&A(n)& \quad \text{ if } p\not | \: (m+2),\\
&A(n)_{\neq x^{\tau}}& \quad \text{ if } p\: | \:(m+2),\\
\end{sis}$$
where $A(n)_{\neq x^{\tau}}:=\{x^a\in A(n)\: |\: x^a\neq x^{\tau}:=x_1^{p-1}\cdots 
x_n^{p-1}\}$. 
Therefore $K(n)$ has dimension $p^n$ if $p\not | \: (m+2)$ and $p^n-1$ if  $p\: | \:(m+2)$.
\end{pro}
\begin{proof}
See \cite[Chapter 4.5]{FS}.
\end{proof}

Kostrikin and Shafarevich (in the above mentioned paper \cite{KS}) conjectured 
that a restricted simple Lie algebras (that is a restricted algebras without proper ideals)
over an algebraically closed field of 
characteristic $p>5$ is either of classical or Cartan type. The {\bf Kostrikin-Shafarevich
conjecture} was proved by Block-Wilson (see \cite{BW1} and \cite{BW2})
for $p>7$, building upon the work of Kostrikin-Shafarevich
(\cite{KS} and \cite{KS2}), Kac (\cite{KAC1} and \cite{KAC3}), Wilson (\cite{WIL2}) 
and Weisfailer (\cite{WEI}).

Recently, Premet and Strade (see \cite{PS1}, \cite{PS2}, \cite{PS3}, \cite{PS4}) 
proved the Kostrikin-Shafarevich conjecture for $p=7$. Moreover they  showed
that for $p=5$ there is only one exception, the Melikian algebra (\cite{MEL}), whose 
definition is given below.
  
\begin{defi}
Let $p={\rm char}(F)=5$. Let $\ti{W(2)}$ be a copy of $W(2)$
and for an element $D\in W(2)$ we indicate with $\widetilde{D}$ the corresponding element
inside $\widetilde{W(2)}$.
The Melikian algebra $M$ is defined as
$$M=A(2)\oplus W(2)\oplus \widetilde{W(2)},$$
with Lie bracket defined by the following rules (for all $D, E\in W(2)$ and 
$f, g\in A(2)$):
$$\begin{sis}
&[D,\widetilde{E}]:=\widetilde{[D,E]}+2\,{\rm div}(D)\widetilde{E},\\
&[D,f]:=D(f)-2\, {\rm div}(D)f,\\
&[f_1\widetilde{D_1}+f_2\widetilde{D_2},g_1\widetilde{D_1}+g_2\widetilde{D_2}]:=
f_1g_2-f_2g_1,\\
&[f,\widetilde{E}]:=f E,\\
&[f,g]:=2\,(gD_2(f)-fD_2(g))\widetilde{D_1}+2\,(fD_1(g)-gD_1(f))\widetilde{D_2},\\
\end{sis}$$
where ${\rm div}(f_1D_1+f_2D_2):=D_1(f_1)+D_2(f_2)\in A(2)$.
\end{defi}

In characteristic $p=2, 3$, there are many exceptional restricted simple Lie algebras 
(see \cite[page 209]{STR}) and the classification seems still far away.

\section{Infinitesimal Deformations}

An infinitesimal deformation of a Lie algebra 
$L$ over a field $F$ is a Lie algebra $L'$ over $F[\epsilon]/(\epsilon^2)$ such that
$L'\times_{F[\epsilon]/(\epsilon^2)}F\cong L$. Explicitly, $L'=L+\epsilon L$
with Lie bracket $[-,-]'$ defined by 
(for any two elements $X,Y\in L\subset L'$): 
$$[X,Y]'=[X,Y]+\epsilon f(X,Y),$$ 
where $[-.-]$ is the Lie bracket of $L$ and $f(-,-)$ is an $2$-alternating function 
from $L$ to $L$, considered a module over itself via adjoint representation. 
The Jacobi identity for $[-,-]'$ forces $f$ to be a cocycle and 
moreover one can check that two cocycles differing by a coboundary define isomorphic 
Lie algebras. Therefore the infinitesimal deformations of a
Lie algebra $L$ are parametrized by
the second cohomology $H^2(L,L)$ of the Lie algebra with values in the adjoint 
representation (see \cite{GER1} for a rigorous treatment).

It is a classical result that simple Lie algebras in characteristic zero 
are rigid. We want to give a sketch of the proof of the following Theorem (see 
\cite{HiSt} for details).

\begin{theo}
Let $L$ be a simple Lie algebra over a field $F$ of characteristic $0$.
Then, for every $i\geq 0$, we have that 
$$H^i(L,L)=0.$$ 
\end{theo}
\begin{proof}[Sketch of the Proof]
Since the Killing form $\beta(x,y)={\rm tr}({\rm ad}(x){\rm ad}(y))$ 
is non-degene\-rate (by Cartan's criterion), we can choose 
two bases $\{e_i\}$ and $\{e_i'\}$ of $L$ such that 
$\beta(e_i,e_j')=\delta_{ij}$. 
Consider the Casimir element
$C:=\sum_i e_i\otimes e_i'$ inside the enveloping algebra $\mathfrak{U}(L)$.
One can check that:
\begin{enumerate}
\item $C$ belongs to the center of the enveloping algebra
and therefore it induces an $L$-homomorphism $C:L\to L$, where $L$ is a consider 
a module over itself via adjoint action. Moreover since 
${\rm tr}_L(C)={\rm dim}(L)\neq 0$, $C$ is non-zero and hence is an isomorphism
by the simplicity of $L$.
\item The map induced by $C$ on the exact complex 
$\{\mathfrak{U}(L)\otimes_k \bigwedge^n L\}_n\to F$ is homotopic to $0$.
\end{enumerate}
Therefore the induced map on cohomology $C_*:H^*(L,L)\to H^*(L,L)$ is an 
isomorphism by $(1)$ and the zero map by $(2)$, which implies that $H^*(L,L)=0.$
\end{proof}

The above proof uses the
non-degeneracy of the Killing form and the non-vanishing of the
trace of the Casimir element, which is equal to the dimension of the
Lie algebra. Therefore the same proof works also for the restricted simple
Lie algebras of {\bf classical type} over a field of
characteristic not dividing the determinant of the Killing form and
the dimension of the Lie algebra. Actually Rudakov (see \cite{RUD})
showed that such Lie algebras are rigid if the characteristic of the
base field is greater or equal to $5$ while in characteristic $2$
and $3$ there are non-rigid classical Lie algebras (see \cite{Che1},
\cite{Che2}, \cite{Che3}).

It was already observed by  Kostrikin and D\v zumadildaev (\cite{DK},
\cite{DZ1}, \cite{DZ2} and \cite{DZ3}) that Witt-Jacobson Lie algebras  
admit infinitesimal deformations.
More precisely: in \cite{DK} the authors compute the infinitesimal
deformations of the Jacobson-Witt algebras of rank $1$, while in
\cite[Theorem 4]{DZ1}, \cite{DZ2} and \cite{DZ3} the author describes the 
infinitesimal deformations of the Jacobson-Witt algebras of any rank but without a
detailed proof.

In the papers \cite{VIV1}, \cite{VIV2} and \cite{VIV3}, we computed the infinitesimal 
deformations of the restricted simple Lie algebras of {\bf Cartan type} in characteristic 
$p\geq 5$, showing in particular that they are non-rigid.
Before stating the results, we need to recall the definition of the Squaring operators
(\cite{GER1}).
The Squaring of a derivation $D:L\to L$ is the $2$-cochain defined, for any $x, y\in L$
as it follows  
\begin{equation}\label{Square}
\Sq(D)(x,y)=\sum_{i=1}^{p-1}\frac{[D^i(x),D^{p-i}(y)]}{i!(p-i)!},
\end{equation}
where $D^i$ is the $i$-iteration of $D$. Using the Jacobi identity, it is 
straightforward to check that $\Sq(D)$ is a $2$-cocycle and therefore it defines
a class in the cohomology group $H^2(L,L)$, which we will continue to call 
$\Sq(D)$ (by abuse of notation). Moreover for an element $\gamma\in L$, we define
$\Sq(\gamma):=\Sq({\rm ad}(\gamma))$.

\begin{theo}\label{W-finaltheorem}
We have that
$$H^2(W(n),W(n))=\bigoplus_{i=1}^n \langle {\rm Sq}(D_i)\rangle_F.$$
\end{theo}

\begin{theo}\label{S-finaltheorem}
We have that
$$H^2(S(n),S(n))=\bigoplus_{i=1}^n  \langle {\rm Sq}(D_i)\rangle_F   \bigoplus
\langle\Theta \rangle_F,$$
where $\Theta$ is defined by $\Theta(D_i,D_j)=D_{ij}(x^{\tau})$.
\end{theo}

\begin{theo}\label{K-finaltheorem}
Let $n=2m+1\geq 3$. Then we have that
$$H^2(K(n),K(n))=\bigoplus_{i=1}^{2m} \langle {\rm Sq}(D_K(x_i)) \rangle_F \bigoplus 
\langle {\rm Sq}(D_K(1))\rangle_F.$$
\end{theo}

Before stating the next theorem, we need some notations about $n$-tuples of
natural numbers. We consider the order relation inside $\N^n$ given by 
$a=(a_1,\cdots, a_n)\leq b=(b_1,\cdots,b_n)$ if $a_i\leq b_i$ for every $i=1,\cdots, n$.  
We define the degree of $a\in \N^n$ as $|a|=\sum_{i=1}^n a_i$ and the 
factorial as $a!=\prod_{i=1}^n a_i!$. For two multindex $a, b\in \N^n$ such that 
$b\leq a$, we set 
$\binom{a}{b}:=\prod_{i=1}^n \binom{a_i}{b_i}=\frac{a!}{b!(a-b)!}.$
For every integer
$j\in \{1,\cdots,n\}$ we call $\ep_j$ the $n$-tuple having $1$ at the $j$-th entry 
and $0$ outside. We denote with $\sigma$ the multindex $(p-1, \cdots, p-1)$.

Assuming now that $n=2m$, we define the sign 
$\sigma(j)$ and the conjugate $j'$ of $1\leq j\leq 2m$ as follows:
$$\sigma(j)=\begin{sis}
1& \quad \text{ if }1\leq j\leq m,\\
-1& \quad \text{ if } m< j\leq 2m,\\
\end{sis}
\hspace{0,5cm} \text{ and } \hspace{0,5cm}
j'=\begin{sis}
j+m& \quad \text{ if } 1\leq j\leq m,\\
j-m& \quad \text{ if } m<j\leq 2m.
\end{sis}$$
Given a multindex $a=(a_1,\cdots,a_{2m})\in \N^{2m}$, we define the sign of $a$ as
$\sigma(a)=\prod\sigma(i)^{a_i}$ and the conjugate of $a$ as the multindex 
$\hat{a}$ such that $\hat{a}_i=a_{i'}$ for every $1\leq i\leq 2m$.

 \begin{theo}\label{H-finaltheorem}
Let $n=2m\geq 2$. Then if $n\geq 4$ we have that
$$H^2(H(n),H(n))=\bigoplus_{i=1}^n \langle {\rm Sq}(D_H(x_i))\rangle_F 
\bigoplus_{\stackrel{i<j}{j\neq i'}} \langle \Pi_{ij} \rangle_F \bigoplus_{i=1}^m 
\langle \Pi_i \rangle_F \bigoplus \langle \Phi\rangle_F,$$
where the above cocycles are defined (and vanish outside) by 
$$\begin{sis}
& \Pi_{ij}(D_H(x^a), D_H(x^b))=D_H(x_{i'}^{p-1}x_{j'}^{p-1}[D_i(x^a)D_j(x^b)-
D_i(x^b)D_j(x^a)]),\\
& \Pi_i(D_H(x_ix^a),D_H(x_{i'}x^b))=D_H(x^{a+b+(p-1)\ep_i+(p-1)\ep_{i'}}), \\
& \Pi_i(D_H(x_k),D_H(x^{\sigma-(p-1)\ep_i-(p-1)\ep_{i'}}))=-\sigma(k)D_H(x^{\sigma-\ep_{k'}}) 
\: \text{ for } 1\leq k\leq n,\\
& \Phi(D_H(x^a),D_H(x^b))=\sum_{\stackrel{\delta\leq a, \widehat{b}}{|\delta|=3}}
\binom{a}{\delta}\binom{b}{\widehat{\delta}}\sigma(\delta)\:\delta!\:
D_H(x^{a+\widehat{b}-\delta -\widehat{\delta}}).\\
\end{sis}$$
If $n=2$ then we have that 
$$H^2(H(2),H(2))=\bigoplus_{i=1}^2 \langle {\rm Sq}(D_H(x_i))\rangle_F 
\bigoplus \langle \Phi\rangle_F.$$

\end{theo}

\begin{theo}\label{M-finaltheorem}
We have that
$$H^2(M,M)=\langle {\rm Sq}(1)\rangle_F \bigoplus_{i=1}^2 \langle {\rm Sq}(D_i)
\rangle_F \bigoplus_{i=1}^2 {\rm Sq}(\ti{D_i})\rangle_F.$$
\end{theo}

\section{Open Problems}

{\bf Simple Lie algebras} (not necessarily restricted) over an algebraically closed field 
$F$ of characteristic $p\neq 2,3$ have been classified by Strade and Wilson for $p>7$
(see  \cite{SW}, \cite{STR1}, \cite{STR2}, \cite{STR3}, \cite{STR4}, \cite{STR5}, 
\cite{STR6}) and by Premet-Strade for $p=5, 7$ (see \cite{PS1}, \cite{PS2}, \cite{PS3},
 \cite{PS4}). The classification says that for $p\geq 7$ a simple Lie algebra is
of classical type (and hence restricted) or of generalized Cartan type. 
Those latter are generalizations of the Lie algebras of Cartan type, obtained 
by considering higher truncations of divided power algebras (not just $p$-truncated 
polynomial algebras) and by considering only the subalgebra of (the so called) 
special derivations (see \cite{FS} or \cite{STR} for the precise definitions).
Again in characteristic $p=5$, the only exception is represented by the 
generalized Melikian algebras. Therefore an interesting problem would be 
the following:

\begin{prob}
Compute the infinitesimal deformations of the simple Lie algebras.
\end{prob}

Note that there is an important distinction between {\bf restricted simple} Lie algebras
and {\bf simple restricted} Lie algebras.
The former algebras are the restricted Lie algebras which do not have any nonzero 
proper ideal, while the second ones are the restricted Lie algebras 
which do not have any nonzero proper restricted ideal (or $p$-ideal), 
that is an ideal closed under the $p$-map. 
Clearly every restricted simple Lie algebra is a simple restricted Lie algebra, but a simple 
restricted Lie algebra need not be a simple Lie algebras. Indeed we have the following

\begin{pro}
There is a bijection
$$\{\text{Simple restricted Lie algebras}\}\longleftrightarrow 
\{\text{Simple Lie algebras}\}.$$
Explicitly to a simple restricted Lie algebra $(L,[p])$ we associates its derived algebra
$[L,L]$. Conversely to a simple Lie algebra $M$ we associate the restricted subalgebra 
$M^{[p]}$ of 
${\rm Der}_F(M)$ generated by ${\rm ad}(M)$ (which is called the universal $p$-envelope of 
$M$).
\end{pro}
\begin{proof}
We have to prove that the above maps are well-defined and are inverse one of the other.

$\bullet$ Consider a simple restricted Lie algebra $(L,[p])$. The derived subalgebra 
$[L,L]\lhd L$ is a non-zero ideal (since $L$ can not be abelian) and therefore 
$[L,L]_p=L$, where $[L,L]_p$ denotes the $p$-closure of $[L,L]$ inside $L$.

Take a non-zero ideal $0\neq I\lhd [L,L]$. 
Since $[L,L]_p=L$, we deduce from \cite[Chapter 2, Prop. 1.3]{FS} that $I$ is also an 
ideal of $L$ and therefore $I_p=L$ by restricted simplicity of $(L,[p])$.
From loc. cit., it follows also that $[L,L]=[I_p,I_p]=[I,I]\subset I$ from which 
we deduce that $I=L$. Therefore $[L,L]$ is simple. 

Since ${\rm ad}:L\to {\rm Der}_F(L)$ is injective and $[L,L]_p=L$, it follows 
by loc. cit. that ${\rm ad}:L\to {\rm Der}_F([L,L])$ is injective. Therefore we have that
$[L,L] \subset L \subset {\rm Der}_F([L,L])$ and hence $[L,L]^{[p]}=[L,L]_p=L$.

$\bullet$ Conversely, start with a simple Lie algebra $M$ and consider its universal
$p$-envelop $M<M^{[p]}<{\rm Der}_F(M)$. 

Take any restricted ideal $I\lhd_p M^{[p]}$. By loc. cit., 
we deduce $[I,M^{[p]}]\subset I\cap [M^{[p]}, M^{[p]}]=I\cap [M,M]=I\cap M\lhd M$. 
Therefore, by the simplicity of $M$, either $I\cap M=M$ or $I\cap M=0$. 
In the first case, we have that $M\subset I$ and therefore $M^{[p]}=I$.
in the second case, we have that $[I,M^{[p]}]=0$ and therefore $I=0$
because $M^{[p]}$ has trivial center. We conclude that $M^{[p]}$ is simple restricted.

Moreover, by loc. cit., we have that $[M^{[p]},M^{[p]}]=[M,M]=M$.

\end{proof}

Therefore the preceding classifications of simple Lie algebras (for $p\neq 2,3$)
give a classification of simple restricted Lie algebras.

\begin{prob}\label{prob2}
Compute the infinitesimal deformations of the simple restricted Lie algebras.
\end{prob}

There is an important connection between simple restricted Lie algebras and 
{\bf simple finite group schemes}.

\begin{pro}
Over an algebraically closed field $F$ of characteristic $p>0$, a simple finite group 
scheme is either the constant group scheme associated to a simple finite group 
or it is the finite group scheme of height $1$ associated to a simple restricted Lie algebra. 
\end{pro}
\begin{proof}
Let $G$ be a simple finite group scheme. The kernel of the Frobenius map $F:G\to G^{(p)}$
is a normal subgroup and therefore, by the simplicity of $G$,
we have that either ${\rm Ker}(F)=0$ or ${\rm Ker}(F)=G$. In the first case, the group $G$ is 
constant (since $F=\overline{F}$), and therefore it corresponds to an (abstract) simple
finite group. In the second case, the group $G$ is of height $1$ and therefore the
result follows from Proposition \ref{height-1}.
\end{proof}

The following problem seems very interesting.

\begin{prob}\label{prob3}
Compute the infinitesimal deformations of the simple finite group schemes.
\end{prob}

Since constant finite group schemes (or more generally \'etale group schemes) are rigid,  
one can restrict to the simple finite group schemes of height $1$ associated to 
the simple restricted Lie algebras. Moreover, if $(L,[p])$ is the simple restricted Lie
algebra corresponding to the simple finite group scheme $G$, then the infinitesimal 
deformations of $G$ correspond to restricted infinitesimal deformations of 
$(L,[p])$, that are infinitesimal deformations that admit a restricted structure. 
These are parametrized by the second restricted cohomology group $H^2_*(L,L)$
(defined in \cite{HOC}). Therefore the above Problem \ref{prob3} is equivalent 
to the following: 

\begin{prob}\label{prob4}
Compute the restricted infinitesimal deformations of the simple restricted Lie algebras.
\end{prob}

The above Problem \ref{prob4} is closely related to Problem \ref{prob2} because
of the following spectral sequence relating the restricted cohomology to the ordinary 
one (see \cite{FAR}):
$$ E_2^{p,q}={\rm Hom}_{{\rm Frob}}\left(\bigwedge^q L, H_*^p(L,L)\right)\Rightarrow 
H^{p+q}(L,L),$$
where ${\rm Hom}_{{\rm Frob}}$ denote the homomorphisms that are semilinear with respect
to the Frobenius.

\begin{ackno}
The results presented here constitute my PHD thesis. I thank my supervisor prof.
R. Schoof for useful advices and constant encouragement.
\end{ackno}


\begin{thebibliography}{INTRO}

\addcontentsline{toc}{section}{References}

\bibitem[AF54]{AF} A. A. Albert and M. S. Frank: \emph{Simple Lie algebras of 
characteristic $p$}. Rend. Sem. Mat. Univ. Politec. Torino {\bf 14} (1954-1955), 117--139.

\bibitem[BAE02]{BAE} J. Baez: \emph{The Octonions}. Bull. Amer. Math. Soc. 39 (2002), 
145-205.

\bibitem[BLO62]{BLO} R. E. Block: \emph{Trace forms on Lie algebras}. Canad. J. Math.
{\bf 14} (1962), 553--564.

\bibitem[BW84]{BW1} R. E. Block and R. L. Wilson: \emph{The restricted simple Lie algebras
are of classical or Cartan type}. Proc. Nat. Aca. Sci. USA {\bf 81} (1984), 5271--5274.

\bibitem[BW88]{BW2} R.E. Block and R.L. Wilson: \emph{Classification of the restricted simple
Lie algebras}.  J. Algebra  114  (1988), 115--259.

\bibitem[CAR09]{CAR} E. Cartan: \emph{Les groupes de transformations continus, 
 infinis, simples}. Ann. Sci. \'Ecole Nor. Sup. {\bf 26} (1909), 93--161.

\bibitem[CHE05]{Che1} N. G. Chebochko: \emph{Deformations of classical Lie algebras with a
 homogeneous root system in characteristic two. I} (Russian).
Mat. Sb.  196  (2005),  no. 9, 125--156. English translation:  Sb. Math.  196  (2005),
no. 9-10, 1371--1402.

\bibitem[CK00]{Che2} N. G. Chebochko and M. I. Kuznetsov: \emph{Deformations of classical Lie
algebras} (Russian). Mat. Sb. 191 (2000), no. 8, 69--88.
English translation: Sb. Math. 191 (2000), no. 7-8, 1171--1190.

\bibitem[CKK00]{Che3} N. G. Chebochko, S. A. Kirillov and M. I. Kuznetsov:
\emph{Deformations of a Lie algebra of type $G\sb 2$ of characteristic three (Russian)}.
Izv. Vyssh. Uchebn. Zaved. Mat. 2000, no. 3, 33--38.
English translation: Russian Math. (Iz. VUZ) 44 (2000), no. 3, 31--36.





\bibitem[DG70]{DG} M. Demazure and P. Gabriel: \emph{Groupes alg\'ebriques}.
Tome I: G\'eom\'etrie alg\'ebrique, g\'en\'eralit\'es, groupes commutatifs (French).
Avec un appendice  Corps de classes local par Michiel Hazewinkel.
Masson and Cie, Editeur, Paris; North-Holland Publishing Co., Amsterdam, 1970.

\bibitem[DK78]{DK} A. S. D\v zumadildaev and A. I. Kostrikin:
\emph{Deformations of the Lie algebra $W\sb{1}(m)$} (Russian).
Algebra, number theory and their applications.  Trudy Mat. Inst. Steklov.  148  (1978),
141--155, 275.

\bibitem[DZU80]{DZ1} A. S. D\v zumadildaev: \emph{Deformations of general Lie algebras
of Cartan type} (Russian).  Dokl. Akad. Nauk SSSR  251  (1980), no. 6, 1289--1292.
English translation: Soviet Math. Dokl. 21 (1980), no. 2, 605--609.

\bibitem[DZU81]{DZ2} A. S. D\v zumadildaev: \emph{Relative cohomology and deformations
of the Lie algebras of Cartan types} (Russian).  Dokl. Akad. Nauk SSSR  257  (1981), no. 5,
1044--1048.
English translation: Soviet Math. Dokl. 23 (1981), no. 2, 398--402.

\bibitem[DZU89]{DZ3} A. S: D\v zumadildaev: \emph{Deformations of the Lie algebras
$W\sb n(m)$} (Russian).  Mat. Sb.  180  (1989),  no. 2, 168--186.
English translation: Math. USSR-Sb.  66  (1990),  no. 1, 169--187.


\bibitem[FAR91]{FAR} R. Farnsteiner: \emph{Cohomology groups of reduced enveloping algebras}.  
Math. Zeit.  206  (1991), 103--117. 


\bibitem[FS88]{FS} R. Farnsteiner and H. Strade: \emph{Modular Lie algebras
and their representation}. Monographs and textbooks in pure and applied mathematics,
vol. 116. Dekker, New York, 1988.

\bibitem[FRA54]{FRA1} M. S. Frank: \emph{A new class of simple Lie algebras}.
Proc. Nat. Acad. Sci. USA {\bf 40} (1954), 713--719. 

\bibitem[FRA64]{FRA2} M. S. Frank: \emph{Two new classes of simple Lie algebras}.
Trans. Amer. Math. Soc. {\bf 112} (1964), 456--482. 

\bibitem[GER64]{GER1} M. Gerstenhaber: \emph{On the deformation of rings and algebras}.
Ann. of Math. 79  (1964), 59--103.





\bibitem[HS97]{HiSt} P. J. Hilton and U. A. Stammbach: \emph{A course in homological 
algebra}. Second edition. Graduate Texts in Mathematics, 4. Springer-Verlag, New York, 1997.

\bibitem[HOC54]{HOC} G. Hochschild: \emph{Cohomology of restricted Lie algebras}.
Amer. J. Math. 76 (1954), 555--580.


\bibitem[JAC37]{JAC2} N. Jacobson: \emph{Abstract derivations and Lie algebras}.
Trans. Amer. Math. Soc. {\bf 42} (1937), 206--224.

\bibitem[JAC80]{JACALG} N. Jacobson: \emph{Basic algebra II}. 
W. H. Freeman and Co., San Francisco, Calif., 1980.


\bibitem[JAC43]{JAC1} N. Jacobson: \emph{Classes of restricted Lie algebras of 
characteristic $p$, II}. Duke Math. Journal {\bf 10} (1943), 107-121.  

\bibitem[JAC71]{JACEXE} N. Jacobson: \emph{Exceptional Lie algebras}. 
Lecture Notes in Pure and Applied Mathematics {\bf 1}, Marcel Dekker, New York 1971.


\bibitem[JAC62]{JACLIE} N. Jacobson: \emph{Lie algebras}. 
Interscience Tracts in Pure and Applied Mathematics No. 10, New York-London, 1962.


\bibitem[KAC70]{KAC1} V. G. Kac: \emph{The classification of the simple Lie algebras
over a field with nonzero characteristic} (Russian). Izv. Akad. Nauk SSSR Ser. Mat. {\bf 34}
(1970), 385--408. English translation: Math. USSR-Izv. {\bf 4} (1970), 391--413.


\bibitem[KAC74]{KAC3} V. G. Kac: \emph{Description of the filtered Lie algebras with which 
graded Lie algebras of Cartan type are associated} (Russian).  
Izv. Akad. Nauk SSSR Ser. Mat. {\bf 38} (1974), 800--834. 
English translation: Math. USSR-Izv. {\bf 8} (1974), 801--835.

\bibitem[KAP71]{KAP} I. Kaplansky: \emph{Lie algebras and locally compact groups}.
Univ. of Chicago Press, Chicago, 1971.

\bibitem[KS66]{KS} A. I. Kostrikin and I. R. Shafarevich: \emph{Cartan's pseudogroups
and the $p$-algebras of Lie} (Russian).  Dokl. Akad. Nauk SSSR  168  (1966), 740--742.
English translation: Soviet Math. Dokl. 7 (1966), 715--718.

\bibitem[KS69]{KS2} A. I. Kostrikin and I. R. Shafarevich: \emph{Graded Lie algebras of
finite characteristic} (Russian). Izv. Akad. Nauk SSSR Ser. Math. {\bf 33} (1969), 
251--322. English translation: Math. USSR-Izv. {\bf 3} (1969), 237--304.  

\bibitem[MEL80]{MEL} G. M. Melikian: \emph{Simple Lie algebras of characteristic $5$}
(Russian).  Uspekhi Mat. Nauk  35  (1980), no. 1 (211), 203--204.

\bibitem[PS97]{PS1} A. Premet and H. Strade: \emph{Simple Lie algebras of small characteristic.
I. Sandwich elements}.  J. Algebra  189  (1997), 419--480.

\bibitem[PS99]{PS2} A. Premet and H. Strade: \emph{Simple Lie algebras of small characteristic.
II. Exceptional roots}.  J. Algebra  216  (1999), 190--301.

\bibitem[PS01]{PS3} A. Premet and  H. Strade:\emph{Simple Lie algebras of small characteristic.
III. The toral rank 2 case}.  J. Algebra  242  (2001), 236--337.

\bibitem[PS04]{PS4} A. Premet and H. Strade: \emph{Simple Lie algebras of small characteristic.
 IV. Solvable and classical roots}.  J. Algebra  278  (2004), 766--833.

\bibitem[RUD71]{RUD} A. N. Rudakov: \emph{Deformations of simple Lie algebras} (Russian).
Izv. Akad. Nauk SSSR Ser. Mat. 35 (1971), 1113--1119.


\bibitem[SEL67]{SEL} G. B. Seligman: \emph{Modular Lie algebras}. Ergebnisse der Mathematik
und ihrer Grenzgebiete, Band 40. Springer-Verlag, New York, 1967.


\bibitem[STR89]{STR1} H. Strade: \emph{The Classification of the Simple Modular Lie Algebras: I.
Determination of the two-sections}. Ann. of Math. 130 (1989), 643--677.

\bibitem[STR92]{STR2} H. Strade: \emph{The Classification of the Simple Modular Lie Algebras:
II. The Toral Structure}. J. Algebra 151 (1992), 425--475.

\bibitem[STR91]{STR3} H. Strade: \emph{The Classification of the Simple Modular Lie Algebras:
III. Solution of the Classical Case}. Ann. of Math. 133 (1991), 577--604.

\bibitem[STR93]{STR4} H. Strade: \emph{The Classification of the Simple Modular Lie Algebras:
IV. Determining the Associated Graded Algebra}. Ann. of Math. 138 (1993), 1--59.

\bibitem[STR94]{STR5} H. Strade: \emph{The Classification of the Simple Modular Lie Algebras:
V. Algebras with Hamiltonian Two-sections}. Abh. Math. Sem. Univ. Hamburg 64 (1994), 167--202.

\bibitem[STR98]{STR6} H. Strade: \emph{The classification of the simple modular Lie algebras.
VI. Solving the final case}.  Trans. Amer. Math. Soc.  350  (1998), 2553--2628.

\bibitem[STR04]{STR} H. Strade: \emph{Simple Lie algebras over fields of positive 
characteristic I: Structure theory}. De Gruyter Expositions in Mathematics, 38. 
Walter de Gruyter, Berlin, 2004.

\bibitem[SW91]{SW} H. Strade and R: L. Wilson: \emph{Classification of Simple Lie Algebras over
Algebraically Closed Fields of Prime Characteristic}. Bull. Amer. Math. Soc. 24 (1991), 357--362.

\bibitem[VIV1]{VIV1} F. Vivani: \emph{Deformations of restricted simple Lie algebras I}.
math.RA/0612861.

\bibitem[VIV2]{VIV2} F. Viviani: \emph{Deformations of restricted simple Lie algebras II}.
math.RA/0702499.

\bibitem[VIV3]{VIV3} F. Viviani: \emph{Deformations of the restricted Melikian algebra}.
math.RA/0702594.


\bibitem[WEI78]{WEI} B: Ju. Weisfailer: \emph{On the structure of the minimal ideal of some
graded Lie algebras of characteristic $p>0$}. J. Algebra {\bf 53} (1978), 344--361. 



\bibitem[WIL76]{WIL2} R. L. Wilson: \emph{A structural characterization of the simple
Lie algebras of generalized Cartan type over fields of prime characteristic}.
J. Algebra {\bf 40} (1976), 418--465.



\end{thebibliography}
\end{document}